\documentclass[12pt]{article}
\usepackage{cblhcp-paper}

\def\aatop #1#2{\genfrac{}{}{0pt}{}{#1}{#2}}

\def\aaangle #1{%
  \left\langle \let\\\aatop #1\right\rangle
}

\def\aaparen #1{%
  \left( \let\\\aatop #1\right)%
}
  
\def\aasubset #1{%
  \left\{\let\\\aatop #1\right\}%
}

\def\aasubvar #1{%
  \left\{\let\\\aatop #1\right\}\!\!^\wedge%
}

\def\aacount #1{%
  \left[\let\\\aatop #1\right]%
}

\def\taaangle #1{%
  \begingroup\textstyle\left\langle \let\\\aatop #1\right\rangle\endgroup
}

\def\taaparen #1{%
  \begingroup\textstyle\left(\let\\\aatop #1\right)\endgroup
}

\def\taasubset #1{%
  \begingroup\textstyle\left\{\let\\\aatop #1\right\}\endgroup
}

\def\taasubvar #1{%
  \begingroup\textstyle\left\{\let\\\aatop #1\right\}\!\!^\wedge\endgroup
}

\def\taacount #1{%
  \begingroup\textstyle\left[\let\\\aatop #1\right]\endgroup
}

\def\aasim#1#2{\aacount{\\{#1}{#2}{\sim}}}
\def\taasim#1#2{\taacount{\\{#1}{#2}{\sim}}}

\newcommand\calB{\mathcal{B}}
\newcommand\calG{\mathcal{G}}
\newcommand\calI{\mathcal{I}}
\newcommand\calP{\mathcal{P}}
\newcommand\calR{\mathcal{R}}

\newcommand\calX{\mathcal{X}}
\newcommand\liesl{\mathfrak{sl}}

\newcounter{aaa}
\def\aaa #1%
{%
\par
\goodbreak
\medskip
\refstepcounter{aaa}%
\noindent
\textbf{%
\theaaa\relax
.\,\ %
#1%
}%
\hskip 0.5em\relax
\ignorespaces
}

\makeatletter
\def\theaaa{\textbf{\the\c@aaa}}
\makeatother

\def\cdb#1{%
  \left\langle
    \let\\\aatop
    #1%
  \right\rangle
}

\makeatletter
\newcommand\bbcal[1]{%
  \@bbcal #1%
}

\newcommand\@bbcal[2]{%
  \mathbb{#1}\mathcal{#2}%
}

\makeatother

\begin{document}

\title{Some stumbling first steps towards\\
  linear homology in a nutshell}
\author{Jonathan Fine}
\date{}

\maketitle

\begin{abstract}
  \noindent
In 1985 Bayer and Billera defined a flag vector $f(X)$ for every
convex polytope $X$, and proved some fundamental properties. The flag
vectors $f(X)$ span a graded ring $\calR=\bigoplus_{d\geq0}\calR_d$.
Here $\calR_d$ is the span of the $f(X)$ with $\dim X=d$. It has
dimension the Fibonacci number $F_{d+1}$.

This paper introduces and explores the conjecture, that $\calR$ has a
counting basis $\{e_i\}$. If true then the equation $f(X) = \sum
g_i(X)e_i$ conjecturally provides a formula for the Betti numbers
$g_i(X)$ of a new homology theory. As the $g_i(X)$ are linear
functions of $f(X)$, we call the new theory linear homology.

Further, assuming the conjecture each $g_i$ will have a rank
$r\geq0$. The rank zero part of linear homology will be (middle
perversity) intersection homology. The higher rank $g_i$ measure
successively more complicated singularities. In dimension $d$ we will
have $\dim\calR_d$ linearly independent Betti numbers.

This paper produces a basis $\{e_i\}$ for $\calR$, that is
conjecturally a counting basis.\\
\textbf{Warning.} Conjecture withdrawn in version~2.
\end{abstract}

\kern -20pt
\begingroup
\advance\parskip -5pt
\tableofcontents
\endgroup

\section{Advice to the reader}
\label{sec:advice}

\global\setcounter{aaa}{-1}

\aaa{}
\label{par:zero}
There is a serious error in this paper which prevents it reaching its
goal. In \ref{sec:rank-basis} the `rank basis' and is introduced. It
gives an decomposition of the flag vector ring. I assumed that the
counting basis is consistent with (is a refinement of) this `rank
decomposition'. I've since done calculations which strongly suggest
(or perhaps prove) that this assumption is false.

Thus, a different way of proceeding is needed to obtain the counting
basis. It is no longer possible (or necessary) to rely on the `rank
basis'. In particular, much of the latter part of this paper is a
failed attempt to `square the circle'.

The only changes made in the second version are additions to the title
and abstract, and the addition of this subsection (and of course
changed page breaks). All numbering and cross-references remain the
same, as the first version started at \ref{par:one}.

\aaa{}
\label{par:one}
This paper is the first exposition of a conjectural homology
theory. It is intended for experts, and others who already have a
strong interest in the area.  It is longer than the author would
like. This is because it has dependencies, which are in the process of
being written up \cite{jf:cpr,jf:cplh}. Once done, this paper can be
simplified and shortened.

\aaa{}
The statement [Exercise] has been added to paragraphs that omit a
proof, that a determined reader could provide, based on what is in the
present paper. The exercises are, for example, algebraic and
combinatorial calculations based already stated definitions and
results.

\aaa{}
Terms being defined are \emph{emphasized}.  Sometime definitions are
made incrementally, perhaps first by use or purpose and then
implemention. So the same term might be emphasized several
times. Sometime \emph{emphasis} is used just for emphasis.

\aaa{}
This paper is a journey. In \S\ref{sec:back-and-over} we describe the
starting point, the direction of travel, and the goal. In
\S\ref{sec:product-shadows} and \S\ref{sec:cone-shadows} we construct
the candidate counting basis. This concludes the paper.

\aaa{}
This paper involves four rings, all of which are
isomorphic. Equivalently, we have four different presentations of the
same ring. The rings have a cone operator $C$, and are generated by
cone and product. Thus, the isomorphisms are unique.

\aaa{}
The first ring is the convex polytope flag vector ring. This is
largely due to Bayer and Billera~\cite{bb:gds}. It provides a link to
combinatorial and counting arguments in geometry, particularly convex
polytopes.

\aaa{}
The second ring is the cone product ring. This has a concise algebraic
description, via equations (\ref{eqn:C(U)C(V)}), (\ref{eqn:J(U,V)})
and (\ref{eqn:D}). It allows the use of algebraic methods.

\aaa{}
The third ring, which is only implicitly defined, is the rank basis
ring. It has basis the $\aaangle{I}$ for $I\in\calI$, where $\calI$ is
defined in \ref{par:calI}. Each $I$ is a non-empty sequence of pairs
of counting numbers. It has product as defined in
(\ref{eqn:{00AL}{00BM}}), and cone as in (\ref{eqn:C{n0L_1L'}}) and
(\ref{eqn:C{a{b+1}L}})

\aaa{}
The fourth ring, which is not even implicitly defined, is the
candidate counting ring. It has basis $\aacount{I}$ for $I\in\calI$,
and cone as in \ref{par:[00abL]} and
\S\ref{sec:product-shadows}. Product is defined only by reference to
the rank basis ring.  The author conjectures that this product can be
explicitly defined, in such a way that the coefficients are clearly
counting numbers.

\aaa{}
A word about proofs is helpful here. If a system of constraints has a
single solution, then any two solutions are isomorphic. This follows
from satisfaction of constraints, the isomorphism need not be
explicitly constructed. This can be used to prove, for example,
combinatorial identities.

\aaa{}
We apply this to our four rings. The author has proved that the first
and third rings satisfy the cone product formulas (see \cite{jf:cplh}
and \cite{jf:cpr} respectively). From this it is easy to prove that
they are isomorphic to the cone product ring.

\aaa{}
The author intends to use a similar approach to the fourth ring. Once
the structure coefficients (see \ref{par:def-structure-coeffs}) are
correctly known for the counting ring, all that is left to do is
verify that they satisfy the cone product formula
(\ref{eqn:C(U)C(V)}). This might be hard.

\aaa{}
Thus, any gaps, obscurities or even errors in the present paper are of
less importance, once these structure coefficients have be
obtained. The author suggest that the reader focus on progress towards
the discovery of the counting structure coefficients, and put to one
side any difficulty in the details.

\aaa{}
This paper takes the form of a journey from one ring to the next. The
journey from first to second is done in a single step, namely by
referencing \cite{jf:cplh} (which builds on the ideas in
\cite{jf:MV-IC}).  This is done because no details from this step are
used in the rest of the paper.  Similarly, there is only passing
mention of intersection homology.

\aaa{}
Finally, there is a conjectural fifth ring, which might interest and
help some readers. It is the representation ring of an unknown
algebraic object $\calG$ (see \ref{par:calG-conj} to
\ref{par:calG-insight}). This ring should again provide a connection
to geometry, but this time via the construction of linear homology.

\section{Background and overview}
\label{sec:back-and-over}

\aaa{}
\label{par:def-homology}
By \emph{homology} we mean, in this paper, any systematic way to
construct useful finite dimensional vector spaces. The dimension $n$
of such a vector space is of course a non-negative integer, or in
other words a member of the set $\{0, 1, 2, 3, \ldots\} = \mathbb{N}$
of \emph{counting numbers}.

\aaa{}
\label{par:def-cand-betti-nums}
Put another way, homology is any systematic way to construct useful
integers, that are non-negative because they are the dimension of a
vector space. And sometimes being non-negative is already a useful
property. These integers, until proved non-negative, we will call
\emph{candidate Betti numbers}.

\aaa{}
\label{par:def-linear-homology}
Suppose for each $X\in\calX$ we have a \emph{linear flag vector}
$f(X)$. In this situation we will say that a homology theory is
\emph{linear} if the associated Betti numbers are linear functions of
$f(X)$.

\aaa{}
\label{par:def-flag-vector}
In this paper $\calX$ will be the set $\calP$ of all convex
polytopes. Every $X\in\calX$ has a dimension $\dim X$. Let $\calP_d$
be all $X\in\calP$, with $\dim X = d$.  In 1985 Bayer and Billera
\cite{bb:gds} defined the (linear) \emph{flag vector} $f(X)$, for
$X\in\calP$. Let $\calR_d$ denote the span of $f(X)$, where $X$ ranges
over $\calP_d$.

\aaa{}
\label{par:def-gds-cone}
Bayer and Billera also proved that the sequence $\dim\calR_d$ is the
sequence $1, 1, 2, 3, 5, 8, \ldots$ of Fibonacci numbers. Along the
way they proved that the \emph{cone} (or pyramid) operater
$C:\calP_d\to\calP_{d+1}$ induces, via the equation $C(f(X))=f(C(X))$,
a linear map $C:\calR_d\to\calR_{d+1}$, which we also call \emph{cone}
and denote by $C$.

\aaa{}
\label{par:calR-bilinear-product}
Finally, consider the Cartesian product $X=X\times Y$ of two convex
polytopes. Bayer and Billera also proved, mostly, that the equation
$f(X)f(Y)=f(X\times Y)$ induces a bilinear product and thus a graded
ring structure on $\calR=\bigoplus_{d\geq0}\calR_d$. (They needed a
special case of this result. The argument they used easily extends to
give the general result.)

\aaa{}
\label{par:about-jf:cplh}
All this will be explained in more detail, from the point of view of
linear homology, in the author's preprint \cite{jf:cplh}. In
particular, a geometric construction will explain why the
Bayer-Billera flag vector should rightly be regarded as the
\emph{linear} flag vector (for convex polytopes).

\aaa{}
\label{par:def-fvr}
To summarize, the (convex polytope linear) \emph{flag vector ring}
$\calR=\bigoplus_{d\geq0}\calR_d$ is a graded ring with, of course, a
bilinear product, and a linear cone operator
$C:\calR_d\to\calR_{d+1}$. Further, $\dim\calR_d$ is the sequence of
Fibonacci numbers.

\aaa{}
\label{par:def-cone-simplex}
Let $e_0$ denote the flag vector of a point. We have, of course, that
$e_0=1\in\calR_0$. If $U=C(B)$ for $B,U\in\calR$, we say that $U$ is a
\emph{cone}, or the cone with \emph{base} $B$. We write $e_{n+1} =
C(e_n)$. This defines the \emph{simplices} $e_n$, for $n\geq0$. Every
simplex $e_n$, except $e_0=1$, is also a cone.

\aaa{}
\label{par:def-g_i(X)}
Now suppose that $\{e_i\}$ is a basis for $\calR$, with each $e_i$
lying in some $\calR_d$. The equation
\begin{equation}
  f(X) = \sum g_i(X) e_i
\end{equation}
now defines linear functions $g_i(X)$ of the flag vector $f(X)$. If
$d=\dim X$ then $g_i(X)=0$, unless $e_i\in\calR_d$. Thus, on $\calP_d$
we get a Fibonacci number of linearly independent linear functions.

\aaa{}
\label{par:def-structure-coeffs}
Again using the basis $\{e_i\}$, the equations
\begin{align}
  e_ie_j
  &=\sum\lambda_{ijk}e_k
  \\
  C(e_i)
  &=\sum\mu_{ij}e_j
\end{align}
define the \emph{structure coefficients} $\lambda_{ijk}$ and
$\mu_{ij}$ for $\calR$ (in the $\{e_i\}$ basis).

\aaa{}
\label{par:def-lhn}
The goal of this paper is to define a candidate Betti number formula
for the at present only conjectural linear homology theory (for convex
polytopes). Doing this is what we mean by describing \emph{linear
  homology in a nutshell}. Hence the title of this paper.

\aaa{}
\label{par:def-counting-basis}
To do this we use the structure coefficients. We say that $\{e_i\}$ is
a \emph{counting basis} for $\calR$ if (i)~each $e_i\in\calR_d$ is an
integer linear combination of $f(X)$ for $X\in\calP_d$ and vice
versa. And (ii)~the $\lambda_{ijk}$ are counting numbers (lie in
$\mathbb{N}$). And (iii)~similarly the $\mu_{ij}$ are counting
numbers.

\aaa{}
\label{par:linear-homology-conjecture}
This paper constructs a basis $\{e_i\}$ for $\calR$ that is,
conjecturally, a counting basis. If so, then the author conjectures
that the associated $g_i(X)$ are the Betti numbers for a new way of
constructing useful vector spaces, to be called \emph{linear
  homology}.

\aaa{}
The author hopes that showing $\{e_i\}$ is a counting basis will be a
straightforward calculation, once we have a correct guess for the
structure coefficients. This paper defines some of the candidate
structure coefficients. Product and cone on $\calR$ will determine the
remainder.

\aaa{}
This paper constructs the candidate counting basis incrementally. The
steps are: (i)~The set of $\aaangle{W}$ is the $CD$-basis. (ii)~The
$\aaangle{I}$ give the same basis. (iii)~Here, $I$ is a non-empty
finite sequence of pairs of counting numbers. (iv)~The $\aasubset{I}$
give the rank basis. (v)~The $\aacount{I}$ give the candidate counting
basis.

\aaa{}
We use $e_i$ to denote an element of a general basis. We use
$\aaangle{W}$, $\aaangle{I}$, $\aasubset{I}$, $\aacount{I}$ to
respectively denote an element of the $CD$, $CD$, rank and candidate
counting bases respectively. However, we also use $e_a$ (for
$a\in\mathbb{N}$) to denote the \emph{simplices} in $\calR$. We have
$e_a = \aaangle{C^a} = \aaangle{\\a0} = \aasubset{\\a0} =
\aacount{\\a0}$. We use combinatorial operations on the sequences
$I\in\calI$, and associated algebraic operations.

\aaa{}
\label{par:kunneth-formula}
We can think of the product structure coefficients $\lambda_{ijk}$ as
providing a list of ingredients for constructing the homology of
$Z=X\times Y$ from that of the factors $X$ and $Y$. The $\{e_i\}$
being a counting basis is a strong constraint on this list. K\"unneth,
in his PhD thesis \cite{hk:betti-produkt}, stated and proved correct
the list of ingredients for the product of compact manifolds.

\aaa{}
Similarly, the cone structure coefficients $\mu_{ij}$ should provide a
list of ingredients for the homology of $X=C(B)$ in terms of that of
the base $B$.

\aaa{}
The concept of a ring with a counting basis is not new. Suppose $G$ is
a finite group. Representations of $G$ can be added using direct sum of
vector spaces, and multiplied by using tensor product.

\aaa{}
\label{par:schurs-lemma}
When Schur's Lemma applies, every representation decomposes uniquely
into a direct sum of irreducibles. Thus, from $G$ the representation
ring $\calR_G$ is formed, and the irreducibles provide a counting
basis for this ring. Each structure coefficient is a number, that
counts how many times an irreducible appears in a product.

\aaa{}
\label{par:calG-conj}
The author conjectures that the flag vector ring $\calR$ is isomorphic
to the representation ring $\calR_{\calG}$ for some algebraic object
$\calG$. For this to work, $\calG$ and its representations must have a
cone $C$ operator, and satisfy Schur's Lemma.

\aaa{}
\label{par:calG-think}
Even though $\calG$ is at present unknown, its conjectured existence
provides a way of thinking about $\calR$, which the author finds
helpful. This way of thinking is used in explicitly in
\ref{par:G-e_1e_1}--\ref{par:{0000}{0000}} and implicitly throughout
\S\ref{sec:rank-basis}.

\aaa{}
For non-singular projective algebraic varieties, linear homology
should reduce to the usual homology theory. And in this setting, the
representation ring of the Lie algebra $\liesl_2$ is crucial to the
proof of the hard Lefschetz theorem.

\aaa{}
\label{par:calG-insight}
Thus, $\calG$ might be some sort of amplification of $\liesl_2$, with
a cone operator $C$. Successful construction of
$\calR_{\calG}\cong\calR$ should then provide the desired counting
basis $\{e_i\}$ for $\calR$. However, it may be that it is the
construction of $\{e_i\}$ that provides the insight needed to
construct $\calG$.

\aaa{}
In any case, the goal of this paper is to construct a basis $\{e_i\}$
for $\calR$, that is conjecturally a counting basis. From this follows
conjecturally a formula for linear homology Betti numbers.

\section{The cone product formula}

\aaa{}
In the previous section we saw how linear homology is related to the
problem, of finding a counting basis for the flag vector ring
$\calR$. In this section we provide a purely algebraic description of
$\calR$.

\aaa{}
Recall that $\calR=\bigoplus_{d\geq0}\calR_d$ is a graded ring with a
linear operator $C:\calR_d\to\calR_{d+1}$. Recall also that
$e_0=1\in\calR_0$, and that $e_{n+1}=C(e_n)$ defines the
\emph{simplices} $e_n$, for $n\geq0$.

\aaa{}
We also need the following. Bayer and Billera \cite{bb:gds} showed
that $\dim\calR_d$ gives the Fibonacci sequence $1,1,2,3,5,8$ and so
on. In proving this, they showed that the operators $C$ and product
together generate $\calR$.

\aaa{}
\label{par:cone-product}
Suppose $U,V\in\calR$. In \cite{jf:cpr} the author will prove the
\emph{cone product formula}
\begin{equation}
  \label{eqn:C(U)C(V)}
C(U)C(V) = C(J(U, V)) + DUV
\end{equation}
where
\begin{equation}
  \label{eqn:J(U,V)}
J(U, V) = UC(V) + C(U) V - e_1 UV
\end{equation}
is the \emph{join formula}.

\aaa{}
\label{par:e_1e_1}
Applying (\ref{eqn:J(U,V)}) with $U=V=e_0$ we obtain $J(U,V) = e_0e_1
+ e_1e_0 - e_1e_0e_0 = e_1$. Thus $C(e_0)C(e_0)= C(e_1) + De_0e_0$,
and so $e_1e_1 = e_2 + D$. Thus
\begin{equation}
  \label{eqn:D}
D = e_1 e_1 - e_2
\end{equation}
is forced by (\ref{eqn:C(U)C(V)}) and (\ref{eqn:J(U,V)}).

\aaa{}
We can now define the $CD$-basis for $\calR$. Let $W$ be any finite
word in symbols $C$ and $D$. Use $\emptyset$ to denote the empty
word. The rules (i)~$\aaangle{\emptyset}=e_0$, and
(ii)~$\aaangle{CW}=C(\aaangle{W})$, and
(iii)~$\aaangle{DW}=D\aaangle{W}$ now define $\aaangle{W}$ for any
finite $CD$-word $W$. This is the \emph{$CD$-basis} for $\calR$.

\aaa{}
\label{par:deg-W-d-is-Fibonacci}
We will now show that this is indeed a basis. Writing $\deg C=1$ and
$\deg D=2$ gives a degree $\deg W$ for every $CD$-word. Clearly,
$\aaangle{W}\in\calR_d$ where $d=\deg W$. Further, the number of $W$
with $\deg W=d$ is equal to the Fibonacci number that is
$\dim\calR_d$. So the candidate basis has the correct number of
elements.

\aaa{}
Recall that $C$ and product generate $\calR$. Every element of the
candidate basis, other than $\aaangle{\emptyset}=e_0=1$, is either a
cone on, or $D$ times, a candidate basis element. After removing
leading powers of $D$, any product of candidate basis elements either
has $e_0$ has a factor (and so is trivial), or is a product of cones.

\aaa{}
The cone product formula allows us to write a product of candidate
basis cones as the cone on something plus $D$ times something else. In
both cases, the something is a linear combination of lower degree
candidate basis products. (For $e_1UV$ we use that $e_1=C(e_0)$ is a
cone.)

\aaa{}
Thus, via the cone product formula, any product of candidate basis
elements is a linear combination of candidate basis
elements. Trivially, the same is true for the cone on a candidate
basis element. As cone and product generate $\calR$, the span of the
candidate basis elements is the whole of $\calR$. Counting, as in
\ref{par:deg-W-d-is-Fibonacci}, now shows that the candidate
$CD$-basis is indeed a basis.

\aaa{}
\label{par:alg-def-cpr}
We now have our algebraic description of $\calR$. We can define
$\calR_d$ to be the vector space with basis $\aaangle{W}$, where $\deg
W=d$. Now as usual write $\calR=\bigoplus_{d\geq0}\calR_d$.

\aaa{}
\label{par:alg-def-cpr-ii}
The rule $D\aaangle{W}=\aaangle{DW}$, together with the cone product
formula, define the ring structure on $\calR$. (We also need the
consequence $D=e_1e_1-e_2$ of the cone product formula.) Finally, we
use the rule $C(\aaangle{W})=\aaangle{CW}$ to define the linear
operator $C:\calR_d\to\calR_{d+1}$.

\aaa{}
The arguments in this section, and the existence and properties of the
flag vector ring, show that the rules in \ref{par:alg-def-cpr} and
\ref{par:alg-def-cpr-ii} define a ring $\calR$ (with operator $C$)
that is isomorphic to the flag vector ring.

\aaa{}
The ring defined by the rules in \ref{par:alg-def-cpr} and
\ref{par:alg-def-cpr-ii} we will call the (universal) \emph{cone
  product ring}, and we will also denote it by $\calR$. It is
isomorphic (via $C\mapsto C$) to the flag vector ring of the previous
section, which we will henceforth denote by $\calR_f$.

\aaa{}
Thus, the cone product ring $\calR$ is defined algebraically. The flag
vector ring $\calR_f$ is defined via convex polytopes. By definition,
there is a \emph{flag vector map} $X\mapsto f(X)\in\calR_{f,\dim
  X}$. The arguments in this section, and in particular the cone
product formula, establish an isomorphism $\calR\cong\calR_f$. Thus,
we can think of the flag vector $f(X)$ as lying in the cone product
ring $\calR$.

\section{Simplices}
\label{sec:simplices}

\aaa{}
Most of the rest of this paper is a study of the cone product ring
$\calR$, or in other words the algebraic consequences of the cone
product formula (\ref{eqn:C(U)C(V)}). In this section we focus on
simplices.

\aaa{}
Recall that $e_0=1\in\calR_0$ and $e_{n+1}=C(e_n)\in\calR_{n+1}$
together define the simplices $e_n$ for $n\geq0$. Every simplex,
except $e_0$, is a cone.

\aaa{}
If $e_ae_b$ is a non-trivial product of simplices then $a,b\geq1$ and
so the cone product formula (\ref{eqn:C(U)C(V)}) gives
\begin{equation}
  \label{eqn:e_ae_b=C(J(etc}
  e_a e_b = C(J(e_{a-1}, e_{b-1})) + D e_{a-1} e_{b-1}
\end{equation}
where $J$ is given by the join formula (\ref{eqn:J(U,V)}).

\aaa{}
For $a,b\geq0$ we have
\begin{equation}
  \label{eqn:J(e_ae_b)=...}
  J(e_a,e_b) = e_a e_{b+1} + e_{a+1}e_b - e_1e_ae_b
\end{equation}
and so the formulas (\ref{eqn:e_ae_b=C(J(etc}) and
(\ref{eqn:J(e_ae_b)=...}), for the product and join of two simplices,
are intertwined.

\aaa{}
\label{par:simplex-join-product}
Here is a statement of results. For $a,b\geq0$ we obtain
\begin{equation}
  \label{eqn:J(e_ae_b)=e_a+b+1}
  J(e_a,e_b) = e_{a+b+1}
\end{equation}
while for $a,b\geq1$ we obtain
\begin{equation}
  \label{eqn:e_ae_b=e_{a+b}+etc}
e_a e_b = e_{a+b} + D e_{a-1} e_{b-1}
\end{equation}
which is easily seen to be a consequence of
(\ref{eqn:J(e_ae_b)=e_a+b+1}).

\aaa{}
\label{par:J(e_0,e_b)}
From (\ref{eqn:J(e_ae_b)=...}) we have $J(e_0,e_b)=e_{b+1}$, which
starts the inductive proof of (\ref{eqn:J(e_ae_b)=e_a+b+1}). Note that
for $a=0$, the $e_1UV$ term in (\ref{eqn:J(e_ae_b)=...}) gives a
cancellation, of $e_1e_ae_b$ with $e_{a+1}e_b$.

\aaa{}
Now consider $e_1e_a$ for $a\geq1$. We get $e_1e_a=e_{a+1}+De_{a-1}$,
via (\ref{eqn:e_ae_b=e_{a+b}+etc}). Thus
\begin{equation}
e_1e_ae_b = (e_1e_a)e_b = e_{a+1}e_b + De_{a-1}e_b
\end{equation}
and so
\begin{equation}
  J(e_a,e_b) = e_ae_{b+1} - De_{a-1}e_b
\end{equation}
where here $e_1e_ae_b$ gives a partial cancellation with $e_{a+1}e_b$.

\aaa{}
Assuming (\ref{eqn:e_ae_b=e_{a+b}+etc}) for $e_ae_{b+1}$ we obtain as
required $J(e_a,e_b)=e_{a+b+1}$. Finally, do we have a well founded
induction? In \ref{par:J(e_0,e_b)} we established
(\ref{eqn:J(e_ae_b)=e_a+b+1}) for $a=0$, and hence
(\ref{eqn:e_ae_b=e_{a+b}+etc}) for $a=1$. More generally, given
(\ref{eqn:J(e_ae_b)=e_a+b+1}) for $a=n$ we get
(\ref{eqn:e_ae_b=e_{a+b}+etc}) for $a=n+1$. We have just shown that
(\ref{eqn:e_ae_b=e_{a+b}+etc}) for $a=n$, gives
(\ref{eqn:J(e_ae_b)=e_a+b+1}), also for $a=n$. Thus, the intertwined
induction for both (\ref{eqn:e_ae_b=C(J(etc}) and
(\ref{eqn:J(e_ae_b)=...}) is established.

\aaa{}
Because this proof involved induction only on $a$, it also proves, for
$a\geq1$
\begin{equation}
  \label{eqn:e_aC(U)}
  e_aC(U) = C^{a+1}(U) + D e_{a-1} U
\end{equation}
which is a generalisation of (\ref{eqn:e_ae_b=e_{a+b}+etc}), via
$U=e_{b-1}$. This result is important in \ref{par:e_a[00L]-intro}.

\aaa{}
Here's the proof. As before we have, for $a\geq1$
\begin{equation}
   \label{eqn:e_aC(U)=C(J-etc}
  e_a C(U) = C(J(e_{a-1}, U)) + D e_{a-1} U
\end{equation}
and as $J(e_0,U) = e_0C(U)+e_1U-e_1e_0U = C(U)$ we get
(\ref{eqn:e_aC(U)}) for $a=1$.

\aaa{}
More generally, for $a\geq2$ we have
\begin{equation}
  \label{eqn:J(e_a-1,U)}
J(e_{a-1},U) = e_{a-1} C(U) + e_a U - e_1 e_{a-1} U
\end{equation}
and by induction on $a$ we have
\begin{equation}
  \label{eqn:e_a-1C(U)}
  e_{a-1}C(U) = C^a(U) + D e_{a-1} U
\end{equation}
while from (\ref{eqn:e_ae_b=e_{a+b}+etc}) we have
$e_1e_{a-1}=e_a+De_{a-1}$.

\aaa{}
Substituting these results in (\ref{eqn:J(e_a-1,U)}) we obtain, for
$a\geq1$
\begin{equation}
J(e_{a-1},U) = C^a(U)
\end{equation}
and hence the required formula (\ref{eqn:e_aC(U)}) for $e_aC(U)$. (If
$U=e_b$ then $C^a(U)=e_{a+b}$, as in (\ref{eqn:J(e_ae_b)=e_a+b+1}).)

\aaa{}
\label{par:[ab][cd]}
We now introduce some notation. We will write
\begin{equation}
\taacount{\\ab} = \taaangle{D^bC^a} = D^be_a
\end{equation}
and now the more general formula
\begin{equation}
  \label{eqn:[ab][cd]}
  \taacount{\\ab}
  \taacount{\\cd}
  =
  \taacount{\\{a+c}{b+d}}
  +
  \taacount{\\01}
  \taacount{\\{a-1}b}
  \taacount{\\{c-1}d}
\end{equation}
is a consequence of (\ref{eqn:e_ae_b=e_{a+b}+etc}), provided
$a,c\geq1$.

\aaa{}
The $\aacount{\\ab}$ are part of the candidate counting basis. In
\S\ref{sec:product-shadows} and \S\ref{sec:cone-shadows} we define
candidate counting basis elements
\begin{equation}
  \aacount{
    \\{a_0}{b_0}
    \\{a_1}{b_1}
    \cdots
    \\{a_r}{b_r}
    }
\end{equation}
for $r\geq1$. But first we must define the rank basis.

\section{Rank one examples}
\label{sec:rank-one-examples}

\aaa{}
The $CD$-basis is not a counting basis. The term $-e_1UV$ in the join
formula (\ref{eqn:J(U,V)}) is the cause of difficulty. In the previous
examples, a cancellation took place. In general, it does not. To
resolve this, we have to modify the $CD$-basis.

\aaa{}
The simplest example is $\aaangle{CD}\aaangle{CD}$. We have
\begin{equation}
  J(D, D)
  =
  DC(D) + C(D)D - e_1DD
  =
  2\taaangle{DCD} - \taaangle{DDC} 
\end{equation}
and so
\begin{equation}
  \label{eqn:<CD><CD>}
  \taaangle{CD} \taaangle{CD}
  =
  2 \taaangle{CDCD}
  - \taaangle{CDDC}
  + \taaangle{DDD}
\end{equation}
which contains totally uncancelled the negative $C(e_1UV)$ term
$\aaangle{CDDC}$.

\aaa{}
\label{par:G-e_1e_1}
Recall from \ref{par:calG-conj} that conjecturally each element of the
counting basis is an irreducible representation of $\calG$. We will
now use this, but for guidance only. We don't rely on the conjecture.
The defining equation (see \ref{par:e_1e_1})
\begin{equation}
  e_1e_1 = e_2 + D  
\end{equation}
for $D$ we thus will think of as a decomposition of $e_1e_1$ into
irreducibles.

\aaa{}
\label{par:G-<CD><CD>}
Now consider (\ref{eqn:<CD><CD>}). If $D$ is a representation, then so
are $\aaangle{CD}$ and $\aaangle{CD}\aaangle{CD}$. But the cone
product formula (\ref{eqn:<CD><CD>}) is not yet writing
$\aaangle{CD}\aaangle{CD}$ as a \emph{sum} of
representations. However, if $\aaangle{CDCD}$ or $\aaangle{DDD}$ were
reducible, perhaps something could be done.

\aaa{}
Suppose
\begin{equation}
  \label{eqn:<CDCD>-decomp}
  \taaangle{CDCD}
  =
  (\taaangle{CDCD} - \taaangle{CDDD})
  +
  \taaangle{CDDC}  
\end{equation}
were a representation decomposition of $\aaangle{CDCD}$. In that case
we would have
\begin{equation}
  \taaangle{CD}
  \taaangle{CD}
  =
  2 (\taaangle{CDCD} - \taaangle{CDDD})
  + \taaangle{CDDC}
  + \taaangle{DDD}
\end{equation}
as a representation decomposition of the cone product
$\aaangle{CD}\aaangle{CD}$.

\aaa{}
This is an example of the route we will follow, to cancel the $-e_1UV$
term in $J(U,V)$. Notice that in this example we get positive
cancelling terms from both $UC(V)$ and $C(U)V$, and so the total
effect is to change the sign of $-e_1UV$. In practice, the
justification for following this route is that the subsequent
calculations successfully remove the negative structure coefficients.

\aaa{}
Consider now $\aaangle{CD}$. Suppose, following
(\ref{eqn:<CDCD>-decomp}) we also adopt
\begin{equation}
  \taaangle{CD}
  =
  (\taaangle{CD} - \taaangle{DC})
  + \taaangle{DC}
\end{equation}
as a representation decomposition. (This decomposition also arises
from intersection homology, see \cite{jf:cplh}.)

\aaa{}
In that case, we will also want a representation decomposition of:
\begin{align*}
  (\taaangle{CD}-\taaangle{DC})^2
  \quad&=\quad
  2 \taaangle{CDCD} - \taaangle{CDDC} + \taaangle{DDD}
  \\
  &\qquad
  - 2 \taaangle{DCCD} - 2 \taaangle{DDD}
  \\
  &\qquad
  + \taaangle{DDCC} + \taaangle{DDD}
\end{align*}
(The details of the expansion are left to the reader.)

\aaa{}
Anticipating a notation from the next section, we write:
\begin{align}
  \taasubset{\\00\\00}
  &=
  \taaangle{CD} - \taaangle{DC}
  \\
  \taasubset{\\00\\11}
  &=
  \taaangle{CDDC} - \taaangle{DDCC}
  \\
  \taasubset{\\00\\00\\00}
  &=
  \taaangle{CDCD}
  - \taaangle{DCCD}
  - \taaangle{CDDC}
  + \taaangle{DDCC}
\end{align}

\aaa{}
\label{par:{0000}{0000}}
Using this notation, we obtain
\begin{equation}
  \label{eqn:{0000}{0000}}
  \taasubset{\\00\\00}
  \taasubset{\\00\\00}
  =
  2\taasubset{\\00\\00\\00}
  + \taasubset{\\00\\11}
\end{equation}
as the representation decomposition of
$(\aaangle{CD}-\aaangle{DC})^2$. (The details are left to the reader.)

\aaa{}
Further examples, using
\begin{equation}
  \taasubset{\\00\\ab} = \taaangle{CDD^bC^a} - \taaangle{D^{b+1}C^{a+2}}
\end{equation}
can be done in the same way.

\aaa{}
This is an example in the \emph{rank basis}. See the next section of
the full definition. See \S\ref{sec:rank-product} for the product
formula in the rank basis.

\section{The rank basis}
\label{sec:rank-basis}

\aaa{}
This section introduces the concept of \emph{rank}, and uses it to
define the \emph{rank decomposition}
$\calR=\bigoplus_{r\geq0}\calR_{[r]}$ of the cone product ring. We do
this by explicitly giving the decomposition
\begin{equation}
  \label{eqn:<W>=rank-sum}
\taaangle{W} = \sum_{r\geq0}\taaangle{W}_{[r]}
\end{equation}
of each element $\aaangle{W}$ of the $CD$-basis, where $U_{[r]}$ is
the \emph{rank $r$ part} of $U\in\calR$.

\aaa{}
\label{par:<W>-decomp}
For each $W$ the sum (\ref{eqn:<W>=rank-sum}) is finite, and so has a
top rank component. These components, which can be indexed by
$CD$-words $W$, form the \emph{rank basis} of $\calR$.  This will
become clear in \ref{par:def-rank-leq-2}.

\aaa{}
In the next section we state the product formula in the rank
basis. (This formula will be proved in \cite{jf:cpr}.) The product
structure coefficients will clearly be counting numbers.

\aaa{}
First, we introduce a new way of writing words $W$ in $C$ and $D$. The
basic idea is to use each occurence of $CD$ in $W$ to split $W$ into
pieces. Each piece (less the $CD$) is then of the form $D^bC^a$. The
number of occurences of $CD$ in $W$ is called the \emph{rank} of $W$,
written as $\rank W$.

\aaa{}
If $W$ has rank zero, it is of the form $D^bC^a$. We already write
$\aacount{\\ab}$ for $\aaangle{D^bC^a}$, and the $\aacount{\\ab}$ has
a nice product formula, given by (\ref{eqn:[ab][cd]}). It is the
absence of $CD$ in $W$ that gives this nice formula.
If $W$ has rank one, it is of the form $D^{b_0}C^{a_0} \> CD \>
D^{b_1}C^{a_1}$. Similarly, each rank $r$ word $W$ is given by a
sequence of $(r+1)$ pairs of counting numbers, and vice versa.

\aaa{}
\label{par:calI}
We define the (linear homology) \emph{index set} $\calI$ as
follows. Let $I$ be a sequence of one or more pairs of counting
numbers. We define $\calI$ to be the set of all such sequences. Each
$I\in\calI$ is a sequence of $(r+1)$ pairs (of counting numbers), for
some $r\geq0$. We write $r=\rank I$. We write
$\calI=\bigsqcup_{r\geq0}\calI_{[r]}$ for the resulting disjoint
union.

\aaa{}
From the previous discussion, it follows that each $I\in\calI_{[r]}$
gives a rank $r$ word $W$ in $C$ and $D$, and vice versa. We use this
bijection to define $\aaangle{I}$ for each \emph{index} $I$ in
$\calI$. We define $\deg I$ to be $\deg W$, or in other words
$\aaangle{I}\in\calR_d$ where $d=\deg I$.

\aaa{}
\label{par:def-rank-leq-2}
The following implicitly defines the rank basis, up to rank $2$,
\begin{align}
  \aaangle{\\{a_0}{b_0}}
  &=
  \aasubset{\\{a_0}{b_0}}
  \\
  \aaangle{\\{a_0}{b_0}\\{a_1}{b_1}}
  &=
  \aasubset{\\{a_0}{b_0}\\{a_1}{b_1}}
  +
  \aasubset{\\{a_0 + a_1 + 1}{b_0 + b_1 + 1}}
  \\
  \nonumber
  \aaangle{\\{a_0}{b_0}\\{a_1}{b_1}\\{a_2}{b_2}}
  &=
  \aasubset{\\{a_0}{b_0}\\{a_1}{b_1}\\{a_2}{b_2}}
  \\
  \label{eqn:<>-to-{}-rank-2}
  \nonumber
  &\qquad
  + \aasubset{\\{a_0+a_1+1}{b_0+b_1+1}\\{a_2}{b_2}}
  + \aasubset{\\{a_0}{b_0}\\{a_1+a_2+1}{b_1+b_2+1}}
  \\
  &\qquad
  + \aasubset{\\{a_0+a_1+a_2+2}{b_0+b_1+b_2+2}}
\end{align}
and gives the rank decomposition promised in \ref{par:<W>-decomp}.

\aaa{}
We can simplify this by first writing
\begin{equation}
  I = I_0 I_1 \ldots I_r
\end{equation}
as a sequence $\aaparen{\\{a_i}{b_i}}$ of pairs, and then introducing
a combining operator. We write
\begin{equation}
  \taaparen{\\ab} \# \taaparen{\\cd}
  =
  \taaparen{\\{a+c+1}{b+d+1}}
\end{equation}
and call $\#$ the \emph{mix operator}. (In
$\aaparen{\\ab}\#\aaparen{\\cd}$, the two pairs lose their separate
identity.)

\aaa{}
For example, we can rewrite (\ref{eqn:<>-to-{}-rank-2}) as
\begin{equation}
  \aaangle{L_0L_1L_2}
  =
  \aasubset{L_0L_1L_2}
  + \aasubset{L_0\#L_1L_2}
  + \aasubset{L_0L_1\#L_2}
  + \aasubset{L_0\#L_1\#L_2}
\end{equation}
by using the mix operator.  Like addition, $\#$ is both associative
and commutative.

\aaa{}
This completes the definition of the \emph{rank basis}, whose elements
are $\aasubset{I}$ for $I$ in the index set $\calI$. Clearly, each
$\aaangle{I}$ is the sum of $2^r$ distinct elements of the rank basis,
where $r=\rank I$.

\aaa{}
For clarity, the top row of $I$ gives the exponents of $C$, and the
second row exponents for $D$. Further, the \emph{degree} $\deg I$ of
$I$ is defined by $\aaangle{I}\in\calR_d$. It is the sum of the top
row, plus twice the sum of the second row, plus three times the rank.

\section{The rank basis product formula}
\label{sec:rank-product}

\aaa{}
Products are as nice as possible, in the rank basis. For example
\begin{equation}
  \label{eqn:{ab}{00L}}
  \taasubset{\\ab}\taasubset{\\00L} = \taasubset{\\abL}
\end{equation}
holds for any sequence $L$ of $r\geq0$ pairs of counting
numbers. (This result, and others in this section stated without
proof, will be proved in \cite{jf:cpr}.)

\aaa{}
This is a powerful result. It greatly simplifies product calculations
in the rank basis. It provides a form of \emph{orthogonality}. To help
explain this, we now introduce some notation.

\aaa{}
\label{par:calB}
Every rank $r$ index $I$ in $\calI$ can be written uniquely as $I=AL$,
where $A\in\calI_{[0]}$ is a pair (of counting numbers) and $L$ is a
sequence of $r$ pairs (of counting numbers). We call $A$ the
\emph{head} and $L$ the \emph{body} of the index $I$. We write
$\calB=\bigsqcup_{r\geq0}\calB_{[r]}$ for the set of all possible
\emph{body indexes}.

\aaa{}
Now let $AL$ and $BM$ be two indexes, in head-body form. From
(\ref{eqn:{ab}{00L}}) we have
\begin{equation}
  \taasubset{AL}\taasubset{BM}
  \>\>=\>\>
  \taasubset{A}\taasubset{B}\>\>
  \taasubset{\\00L}\taasubset{\\00M}
\end{equation}
which reduces the general product to computing
$\aasubset{\\00L}\aasubset{\\00M}$. (Because
$\aaangle{\\ab}=\aasubset{\\ab}=\aacount{\\ab}$, we can use
(\ref{eqn:[ab][cd]}) to compute $\aasubset{A}\aasubset{B}$.)

\aaa{}
Recall equation (\ref{eqn:{0000}{0000}})
\begin{equation*}
  \taasubset{\\00\\00}
  \taasubset{\\00\\00}
  =
  2\taasubset{\\00\\00\\00}
  + \taasubset{\\00\\11}
\end{equation*}
from \S\ref{sec:simplices}. This is a special case of a general
property.

\aaa{}
\label{par:def-point-like}
We define the \emph{point-like} elements of $\calR$ as follows. First,
every $\aasubset{\\00L}$ is point-like. Second, any linear combination
of point-like is also point-like. Finally, no other elements are
point-like. In other words, the point-like elements of $\calR$ is the
subspace spanned by the $\aasubset{\\00L}$, for $\calI\in\calB$.

\aaa{}
Also proved in \cite{jf:cpr} is that any product of point-like
elements is also point-like. This is also a powerful result. Equation
(\ref{eqn:{0000}{0000}}), quoted above, is an example of this.

\aaa{}
\label{par:{00L}{00M}}
We will now introduce further notation to help us. By using
\emph{formal sums} the equation
\begin{equation}
  \label{eqn:{00L}{00M}}
  \taasubset{\\00L}\taasubset{\\00M}
  =
  \taasubset{\\00 L*M}
\end{equation}
both (i)~states that the product of point-like is point-like, and
(ii)~defines \emph{merge product} $L*M$, which is a formal sum over
$\calB$.

\aaa{}
Here's how. Suppose $F:S\to V$ is any function, that takes values in a
vector space $V$. We assume no structure on $S$. The equation
\begin{equation}
  F(\lambda_1s_1 + \ldots + \lambda_ns_n)
  =
  \lambda_1F(s_1) + \ldots + \lambda_nF(s_n)
\end{equation}
defines a map, also denoted by $F$, from \emph{formal sums} (of
elements of $S$) to the vector space $V$.

\aaa{}
Note that if $S$ has structure such as addition, care must be taken to
distinguish formal addition from actual addition. For example, if
$F(n)=n^2$ for $n\in\mathbb{N}$, then is $F(2+3)$ equal to $2^2+3^2$,
or is it $(2+3)^2$? Sometime we use $(x)+(y)$ to denote a
\emph{formal} sum.

\aaa{}
We will use this with the vector-valued function
$\aasubset{}:\calI\to\calR$. We return to (\ref{eqn:{00L}{00M}}).
Recall that the product of point-likes is also point-like. Thus we can
write
\begin{equation}
  \taasubset{\\00L}\taasubset{\\00M}
  =
  \lambda_1\taasubset{\\00P_1}
  +
  \ldots
  +
  \lambda_n\taasubset{\\00P_n}
\end{equation}
for suitable $\lambda_i, P_i$.

\aaa{}
Now use these same values to define
\begin{equation}
  P = \lambda_1P_1 + \ldots + \lambda_nP_n
\end{equation}
as a formal sum. By definition,
$\aasubset{\\00L}\aasubset{\\00M}=\aasubset{\\00P}$.

\aaa{}
As promised in \ref{par:{00L}{00M}}, the equation
(\ref{eqn:{00L}{00M}})
\begin{equation*}
  \taasubset{\\00L}\taasubset{\\00M}
  =
  \taasubset{\\00 L*M}
\end{equation*}
now makes two statements. The first is that
$\aasubset{\\00L}\aasubset{\\00M}$ is point-like. The second is that
$L*M$ is defined to be the formal sum, which we previously denoted by
$P$ (for product).

\aaa{}
\label{par:emptyset*M}
We can now state the rank basis point-like product formula. First, the
trivial cases. If $L\in\calB_{[0]}$ then
$\aasubset{\\00L}=\aasubset{\\00}$, and so
$\aasubset{\\00L}\aasubset{\\00M}=\aasubset{\\00M}$. Similarly for
$M\in\calB_{[0]}$.

\aaa{}
\label{par:{00AL}{00BM}}
Any element of $\calB_{[r]}$, for $r\geq1$, can be written uniquely as
$AL$, with $A\in\calB_{[1]}$ and $L\in\calB_{[r-1]}$. In \cite{jf:cpr}
the equation
\begin{equation}
  \label{eqn:{00AL}{00BM}}
  \taasubset{\\00 AL} \taasubset{\\00 BM}
  =
  \taasubset{\\00 A(L*BM)}
  + \taasubset{\\00 B(AL*M)}
  + \taasubset{\\00 (A\#B)(L*M)}
\end{equation}
is proved. This deals with the general case.

\aaa{}
We call $L*M$ the \emph{merge product} of the body indexes $L,
M\in\calB$. The general case equation (\ref{eqn:{00AL}{00BM}}) writes
the merge product $(AL)*(BM)$ as the sum of the \emph{left shuffle}
$A(L*MB)$, the \emph{right shuffle} $B(AL*M)$ and the \emph{mix
  combine} $(A\#B)(L*M)$.

\aaa{}
For clarity, by definition, the equations
\begin{align}
  \taasubset{\\00 L} \taasubset{\\00 BM}
  &=
  \taasubset{\\00 L*BM}
  \\
  \taasubset{\\00 AL} \taasubset{\\00 M}
  &=
  \taasubset{\\00 AL*M}
  \\
  \label{eqn:{00L}{00M}-duplicate}
  \taasubset{\\00 L} \taasubset{\\00 M}
  &=
  \taasubset{\\00 L*M}  
\end{align}
define the quantities $(L*BM)$, $(AL*M)$ and $(L*M)$ in
(\ref{eqn:{00AL}{00BM}}), which thus implicitly contains a recursion.

\aaa{}
For example, if $A,B\in\calB_{[1]}$ then
\begin{equation}
  \label{eqn:{00A}{00B}}
  \taasubset{\\00A}
  \taasubset{\\00B}
  =
  \taasubset{\\00AB}
  + \taasubset{\\00BA}
  + \taasubset{\\00A\#B}
\end{equation}
follows from (\ref{eqn:{00AL}{00BM}}), together with
\ref{par:emptyset*M}. This is a generalisation of
(\ref{eqn:{0000}{0000}}). [Exercise].

\aaa{}
\label{par:rank-r-rank-s}
Suppose $L\in\calB_{[r]}$ and $M\in\calB_{[s]}$, with $r\geq s$. From
(\ref{eqn:{00AL}{00BM}}), say treated as a definition, it follows that
the rank $t$ part $(L*M)_{[t]}$ is zero unless $r\leq t \leq
(r+s)$. Further, within that range the coefficients are counting
numbers, whose sum is non-negative. In fact, the sum depends only on
the pair $r,s$. [Exercise], or see \cite{jf:cpr}.

\aaa{}
Formal sums can also be used for the head part of a product.
Assume $A,B\in\calI_{[0]}$. The equation
\begin{equation}
  \taacount{A}\taacount{B}
  =
  \taacount{\pi(A,B)}
\end{equation}
defines the formal sum $\pi(A,B)$ (which by \S\ref{sec:simplices} is
over $\calI_{[0]}$). As $\aasubset{A}=\aacount{A}$, this is equivalent
to:
\begin{equation}
  \label{eqn:{A}{B}}
  \taasubset{A}\taasubset{B}
  =
  \taasubset{\pi(A,B)}
\end{equation}

\aaa{}
Thus
\begin{equation}
  \taasubset{AL}\aasubset{BM}
  =
  \taasubset{\pi(A,B)\>\> L*M}
\end{equation}
states the rank basis product formula.
In it, the product rule \emph{decomposes orthogonally} into a
\emph{head-product} and a \emph{body-product}, given respectively by
$\pi(A,B)$ and $L*M$.
For clarity, here $A,B\in\calI_{[0]}$ and $L,M\in\calB$. Further,
$\pi(A,B)$ is defined by (\ref{eqn:{A}{B}}), and $L*M$ by
(\ref{eqn:{00AL}{00BM}}).

\section{Rank and cone}
\label{sec:rank-and-cone}

\aaa{}
\label{par:C{0000}-problem}
Although product is counting in the subset basis, cone is not. In the
subset basis, cones often produce negative coefficients, for a simple
reason. For example:
\begin{align}
  \label{eqn:C{0000}}
  \nonumber
  C\taasubset{\\00 \\00}
  &= C\left(\aaangle{CD} - \aaangle{DC}\right)
  = \aaangle{CCD} - \aaangle{CDC}
  \\
  \nonumber
  &=\left(\aaangle{CCD} - \aaangle{DCC}\right)
  -\left(\aaangle{CDC} - \aaangle{DCC}\right)
  \\
  &=
  \taasubset{\\10 \\00} - \taasubset{\\00 \\10}
  =
  \taasubset{\\10}\taasubset{\\00\\00}
  -
  \taasubset{\\00}\taasubset{\\00 \\10}
\end{align}

\aaa{}
\label{par:{1000}=C{0000}+{0010}}
To resolve this particular problem first rewrite (\ref{eqn:C{0000}})
as
\begin{equation}
  \taasubset{\\10}\taasubset{\\00\\00}
  =
  C\taasubset{\\00 \\00}
  +
  \taasubset{\\00}\taasubset{\\00 \\10}
\end{equation}
which then shows, assuming $\calG$, that $\taasubset{\\10 \\00}$ is a
reducible representation. Thus, $\aasubset{\\10\\00}$ is not part of the
counting basis.

\aaa{}
Put another way, the negative coefficient arises because
$\aasubset{\\10 \\00}$ is too large, to be part of a counting
basis. We could fix it by using $C\aasubset{\\00\\00}$ instead of
$\aasubset{\\10\\00}$ as part of the basis.

\aaa{}
In the $CD$ basis, the structure coefficients for $C$ are always
counting numbers, but those for product are not. To solve this we
introduced the rank basis. Now product is always counting, and $C$
sometimes is not. Despite this, we have made progress.

\aaa{}
First, the failure of counting for product in the $CD$ basis is deep
problem, embedded in the $-e_1UV$ term in the join formula. At root,
it is a geometric problem. The failure of cone in the rank basis is a
relatively shallow problem. So going from $CD$ to rank basis solves a
hard problem, at the expense of introducing an easier problem.

\aaa{}
Second, the rank basis introduces the rank decomposition
$\calR=\bigoplus_{r\geq0}\calR_{[r]}$. The $CD$ basis has no similar
feature. (The degree grading $\calR=\bigoplus_{d\geq0}\calR_d$, which
makes $\calR$ a graded ring, is compatible with both bases. In other
words, each basis element $\aaangle{I}$ or $\aasubset{I}$ lies in some
$\calR_d$, where $d=\deg I$).

\aaa{}
Third, the rank decomposition introduces a useful, and perhaps novel,
structure on $\calR$. For example, in \ref{par:rank-r-rank-s} we saw
that rank~$r$ times rank~$s$ is a sum of terms, whose rank lies
between $\max(r,s)$ and $(r+s)$.

\aaa{}
In particular. if $U$ and $V$ are both known up to and including
rank~$r$, then the same is true of the product $UV$ (and trivially
also the sum $U+V$). Thus we can form the \emph{rank~$r$ quotient
  ring} $
\calR/\calR_{[>r]}
\cong
\calR_{[0]}
+ \ldots
+ \calR_{[r]}
$. The author hasn't seen such a structure before.

\aaa{}
Fourth, consider the conjectural $\calG$. If it exists, its
irreducible representations induce a decomposition of $\calR$ into
one-dimensional subspaces. For the rank decomposition to be directly
useful, this $\calG$ decomposition should refine the rank
decomposition.

\aaa{}
In other words, we expect (i)~each irreducible representation $\rho$
to have a rank~$r$, and (ii)~$\calR_{[r]}$ to have as basis the
rank~$r$ representations, and (iii)~when a product of irreducibles is
fully decomposed, the resulting ranks are as in $\calR$, and
(iv)~there is a cone operator $C$ on representations, and (v)~the
representation $C$ operator satisfies the cone product formula.

\aaa{}
Thus, we see that the rank decomposition is conjecturally a step
towards the construction of $\calG$, and it introduces a novel and
interesting mathematical structure.
In the rest of this section, we explore the interaction between $C$
and the rank. First we look at the rank decomposition of the cone
$C(\aasubset{\\00L})$, with point-like base $\aasubset{\\00L}$.

\aaa{}
Consider $C(\aasubset{\\00L})$. If $\rank L=0$, then
$\aasubset{\\00L}=\aasubset{\\00}$ and so $C(\aasubset{\\00L}) =
\aasubset{\\10} = \aasubset{\\10L}$. This takes care of the trivial
case. For the general case of $\rank L\geq0$, we have $L=L_0L'$, where
$L_1=\aatop{a}{b}$ is a pair of counting numbers.

\aaa{}
Now consider $C(\aasubset{\\00L_1L'})$. To evaluate this, we can
unpack into the $CD$ basis (giving an alternating sum), apply $C$, and
then repack. When unpacked the leading rank $CD$ term is
$\aaangle{\\00L_1L'} = C(D\aasubset{L_1L'})$. After applying $C$ we
get $C^2(D\aasubset{L_1L'}) = \aaangle{\\10L_1L'}$. Thus, repacking,
$\aasubset{\\10L_1L'}$ is one of the terms in
$C(\aasubset{\\00L_1L'})$.

\aaa{}
\label{par:C[00L_1L']}
Terms in $\aasubset{\\00L_1L'}$ also arise by applying the mix
operator $\#$ at least once (and alternating the sign). Thus,
$-\aaangle{\\00\#L_1L'}$ appears. But $\aatop00\#\aatop ab =
\aatop{a+1}{b+1}$ and so $C(\aaangle{\\00\#L_1L'})$ is the $\rank r$
element $\aaangle{\\00L_1^{(1)}L'}$, where $L_1^{(1)} =
\aatop{a+1}{b}$. Thus, $\aasubset{\\00L_1^{(0)}}$ is one of the terms
in $C(\aasubset{\\00L_1L'}$.

\aaa{}
The equation
\begin{equation}
  C(\taasubset{\\00 L_1 L'})
  =
  \taasubset{\\10 L_1 L'}
  - \taasubset{\\00 L_1^{(1)} L'}
\end{equation}
now follows, once it is checked that the remaining lower rank terms
cancel. [Exercise]

\aaa{}
\label{par:C{n0L_1L'}}
The same method of proof shows the more general result
\begin{equation}
  \label{eqn:C{n0L_1L'}}
  C(\taasubset{\\n0 L_1 L'})
  =
  \taasubset{\\{n+1}0 L_1 L'}
  - \taasubset{\\00 L_1^{(n)} L'}
\end{equation}
where $L_1^{(n)} = \aatop{a+n}{b}$, when
$L_1=\aatop{a}{b}$. [Exercise]

\aaa{}
\label{par:C{a{b+1}L}}
Similarly, for $a,b\geq0$ the same style of argument shows
\begin{equation}
  \label{eqn:C{a{b+1}L}}
  C(\taasubset{\\{a}{b+1}L})
  =
  \taasubset{\\00\\{a}{b}L}
  + \taasubset{\\{a+1}{b+1}L})
\end{equation}
which completes the calculation of $C$ in the rank basis. [Exercise]

\aaa{}
We now use (\ref{eqn:C{n0L_1L'}}) and (\ref{eqn:C{a{b+1}L}}) to
understand better the interaction between $C$ and the rank
decomposition. First we will give a rank decomposition of the $C$
operator. We then use this decomposition to redefine
\emph{point-like}, and to introduce the related concept
\emph{simplex-like}.

\aaa{}
First, suppose $U\in\calR_{[r]}$. From inspecting
(\ref{eqn:C{n0L_1L'}}) and (\ref{eqn:C{a{b+1}L}}) the equation
\begin{equation}
  C(U) \in \calR_{[r]} \oplus \calR_{[r+1]}
\end{equation}
now follows. Thus, we can write
\begin{equation}
  C = C_{[0]} + C_{[1]}
\end{equation}
where the operation $C_{[i]}$ has rank~$i$ (or in other words is a map
$\calR_{[r]}\to\calR_{[r+1]}$).

\aaa{}
Recall from \ref{par:def-point-like} that $U\in\calR$ is
\emph{point-like} if it is the span of the $\aasubset{\\00L}$. We can
free this concept from its dependence on the rank basis, and also
introduce the related simplex-like. (Key here is that we defined
point-like as a vector subspace of $\calR$, rather than as a property
of basis vectors.)

\aaa{}
First, let's solve the equation $U=C_{[1]}(V)$, for $U$. We start with
$V=\aasubset{AL}$ in the rank basis. If $A=\aatop{a}{0}$ then
(\ref{eqn:C{n0L_1L'}}) applies, and $C(V)=C_{[0]}(V)$, and so
$C_{[1]}(V)=0$. If $A=\aatop{a}{b+1}$ then
$C_{[1]}(V)=\aasubset{\\00\\abL}$.

\aaa{}
This shows that the point-like elements $U$ of $\calR$, as defined in
\ref{par:def-point-like}, are precisely the ones for which the
equation $U=C_{[1]}(V)$ can be solved. In other words, the
\emph{point-like} elements are precisely the \emph{range} of the
operator $C_{[1]}$.

\aaa{}
Now let's solve $C_{[1]}(V)=0$, for $V$. Again, we start with
$V=\aasubset{AL}$ in the rank basis. If $A=\aatop{a}{0}$ then
(\ref{eqn:C{n0L_1L'}}) shows that $C_{[1]}(V)=0$. Further, if
$A=\aatop{a}{b+1}$ then by (\ref{eqn:C{a{b+1}L}}) we have
$C_{[1]}(V)\neq0$. This is almost enough. Some non-trivial linear
combination of not-simplex-like basis vectors might be zero under
$C_{[1]}$.

\aaa{}
To deal with this, note that, if $A=\aatop{a}{b+1}$ then
$C_{[1]}(V)=\aasubset{\\00\\abL}$. Thus, $C_{[1]}$ takes distinct
not-simplex-like elements of the rank basis to distinct elements of
the rank basis. This is enough. It shows that if a linear combination
of not-simplex-like rank basis elements goes to zero under $C_{[1]}$,
then the linear combination is trivial.

\aaa{}
We define the \emph{simplex-like} elements $V$ of $\calR$ to be those
such that $C_{[1]}(V)=0$. In other words, the simplex-like elements
are the \emph{kernel} of the operator $C_{[1]}$. We have shown that
the $\aasubset{\\{a}{0}L}$ provide a basis of the simplex-like
elements.

\aaa{}
In the next section we use equation (\ref{eqn:e_aC(U)})
\begin{equation*}
  e_aC(U) = C^{a+1}(U) + D e_{a-1} U
\end{equation*}
and properties of rank to resume the definition of the counting
basis. In particular, we resolve systematically the problem first
considered in \ref{par:C{0000}-problem}.

\aaa{}
This requires the following result. Recall
$C=C_{[0]}+C_{[1]}$. Consider $C^n=(C_{[0]}+C_{[1]})^n$. Then the
equations
\begin{align}
  \label{eqn:(C^n)_[0]}
  (C^n)_{[0]}
  &=
  (C_{[0]})^n
  \\
  \label{eqn:(C^n)_[1]}  
  (C^n)_{[1]}
  &=
  \sum\nolimits_{i=0}^n
    (C_{[0]})^i
    C_{[1]}
    (C_{[0]})^{n-i}  
\end{align}
give the rank zero and rank one parts of $C^n$. Further, the other
components are zero.

\aaa{}
The proof goes as follows. First, suppose that in the binomial
expansion of $(C_{[0]}+C_{[1]})^n$ we choose $C_{[1]}$ exactly $r$
times. If so, then the resulting product has rank~$r$.  This proves
(\ref{eqn:(C^n)_[0]}) and (\ref{eqn:(C^n)_[1]}).

\aaa{}
Now suppose we choose $C_{[1]}$ at least twice. If so, then the
resulting product contains
\begin{equation}
  C_{[1]} (C_{[0]})^a C_{[1]}
\end{equation}
as a subexpression, for some $a\geq0$. Now, $V=C_{[1]}(U)$ is always
point-like, so $V'=(C_{[0]})^a(V)$ is simplex-like, so $C_{[1]}(V')$
is zero. This proves that the other rank components $(C^n)_{[r]}$ of $C^n$
are zero.

\newpage
\section{Product shadows}
\label{sec:product-shadows}

\aaa{}
The equations
\begin{align}
  \label{eqn:C[0000]}
  C(\taacount{\\00\\00})
  &=
  \taacount{\\10\\00}
  \\
  \label{eqn:[10][0000]}
  \taacount{\\10}
  \taacount{\\00\\00}
  &=
  \taacount{\\10\\00}
  + \taacount{\\00\\10}
\end{align}
are suggested by \ref{par:C{0000}-problem} and
\ref{par:{1000}=C{0000}+{0010}}. They are a special case of equations
and definitions in this section. In (\ref{eqn:[10][0000]}) we call
$\aacount{\\00\\10}$ a \emph{shadow term}. We call
$\aacount{\\10\\00}$ the \emph{leading} term.

\aaa{}
\label{par:[0000]-known}
Assume $\aacount{\\00\\00}$ is known. Equation (\ref{eqn:C[0000]}) now
determines $\aacount{\\10\\00}$. And now (\ref{eqn:[10][0000]})
determines $\aacount{\\00\\10}$. In this section we do something
similar first for $\aacount{\\a0L}$, and then for
$\aacount{\\a{b+1}L}$. The next section determines $\aacount{\\00L}$,
for $\rank L \geq 0$. (If $\rank L=0$,
$\aacount{\\00L}=\aacount{\\00}$, determined by $\aacount{\\00} =
1\in\calR $.)

\aaa{}
The previous section showed that $C=C_{[0]}+C_{[1]}$ is the rank
decomposition of $C$. It defined the \emph{point-wise} subspace of
$\calR$ to be the image of $C_{[1]}$, and the \emph{simplex-like}
subspace to the kernel (or nullspace) of $C_{[1]}$. In the candidate
counting basis, the $\aacount{\\00L}$ will be a basis for the
point-like subspace. Similarly, $\aacount{\\a0L}$ a basis for the
simplex-like.

\aaa{}
Suppose now we are give $\aacount{\\00L}$, for some $L\in\calB$ of
rank $r+1 \geq1$. By assumption (or convention) it is point-like. So
we can solve $\aacount{\\00L}=C_{[1]}(U)$ for $U$. This solution is
unique, up to the addition of a simplex-like $V$. Thus, we can assume
$\rank U=r\geq0$. (For $\rank \aacount{\\00L} = 0$, see
\ref{par:[0000]-known}.)

\aaa{}
\label{par:e_a[00L]-intro}
Now consider $e_a\aacount{\\00L}$, for $a\geq0$. By the rank product
property, this also has rank $r+1$. But $\aacount{\\00L}=C_{[1]}(U)$
and $C=C_{[0]}+C_{[1]}$. So we can apply formula (\ref{eqn:e_aC(U)})
to $e_aC(U)$, to expand $e_a\aacount{\\00L}$.

\aaa{}
We will do this for a typical case, namely $a=2$. We have
\begin{align}
  \label{eqn:pre-e_2[00L]}
  e_2\taacount{\\00L}
  &=
  e_2C(U) - e_2C_{[0]}(U)
  \\
  &=
  C^3(U) + De_1U - e_2C_{[0]}(U)
\end{align}
and now look at the rank $(r+1)$ part.

\aaa{}
\label{par:e_2[00L]}
Using (\ref{eqn:(C^n)_[1]}) and the rank product property we obtain
\begin{equation}
  \label{eqn:e_2[00L]}
  e_2\taacount{\\00L}
  = C_{[0]}C_{[0]}C_{[1]}(U)
  + C_{[0]}C_{[1]}C_{[0]}(U)
  + C_{[1]}C_{[0]}C_{[0]}(U)
\end{equation}
which, we will soon see, determines the shadow terms for
$\aacount{\\20}\aacount{\\00L}$.

\aaa{}
An aside. We obtained (\ref{eqn:e_2[00L]}) by using properties of
rank. This amounts to a massive cancellation of terms in
(\ref{eqn:pre-e_2[00L]}). Perhaps some other massive cancellations
also arise from rank properties, perhaps for an at present unknown
rank.

\aaa{}
To return to (\ref{eqn:e_2[00L]}).  First note that $C_{[1]}$ is
point-like, and so $C_{[0]}C_{[0]}C_{[1]} = C^2C_{[1]}$. In the rank
basis, point-like elements have the form $\aacount{\\00L}$. This leads
to the general rule
\begin{equation}
  \label{eqn:C[a0L]}
  C\taacount{\\a0L} = \taacount{\\{a+1}0L}
\end{equation}
which determines the simplex-like elements, given the point-like.

\aaa{}
Recall that $\aacount{\\00L}=C_{[1]}(U)$ determines $U$ only up to a
simplex-like element $V$. But if $V$ is simplex-like, then so are
$C_{[0]}(V)$, $C_{[0]}C_{[0]}(V)$ and so on. Thus, the
\emph{individual terms} in (\ref{eqn:e_2[00L]}) don't depend on the
choice of $U$. They depend only on $\aacount{\\00L}$.

\aaa{}
Put another way, if $\aacount{\\00L}=C_{[1]}(U)$ then the equation
\begin{equation}
  \taacount{\\00L^{(1)}} = C_{[1]}C_{[0]}(U)
\end{equation}
determines an element $\aacount{\\00L^{(1)}}$ of $\calR$, \emph{which
  does not depend on the choice of $U$}. (This can be though of as a
form of orthogonality.)

\aaa{}
Thus
\begin{equation}
  \taacount{\\00L} \mapsto \taacount{\\00 L^{(1)}}
\end{equation}
defines the \emph{shadow operator} on the point-like subspace of
$\calR$. It preserves rank, and increase degree by one. We use
$L^{(n)}$ to denote the $n$-fold iteration of this operator.

\aaa{}
We call $\aacount{\\00 L^{(n)}}$ the \emph{$n$-th shadow} of
$\aacount{\\00L}$.  Note that $L^{(1)}$ is a formal sum over
$\calB$. This definition agrees with the definition in
\ref{par:C[00L_1L']}.

\aaa{}
We can now rewrite (\ref{eqn:e_2[00L]}) as
\begin{equation}
  \label{eqn:[20][00L]}
  \taacount{\\20}
  \taacount{\\00L}
  =
  \taacount{\\20 L}
  + \taacount{\\10 L^{(1)}}
  + \taacount{\\00 L^{(2)}}
\end{equation}
which, as promised in \ref{par:e_2[00L]}, defines the leading and
shadow terms for $\aacount{\\20}\aacount{\\00L}$. Note that in
(\ref{eqn:[20][00L]}), all terms have the same rank and degree.  The
general case of $\aacount{\\a0}\aacount{\\00L}$ is similar. [Exercise]

\aaa{}
We can rewrite the general case of (\ref{eqn:[20][00L]}) as
\begin{equation}
  \label{eqn:[a0L]}
  \taacount{\\a0L}
  =
  \taacount{\\a0}
  \taacount{\\00L}
  -
  \taacount{\\{a-1}0}
  \taacount{\\00L^{(1)}}
\end{equation}
for $a\geq1$. This writes $\aacount{\\a0L}$ as a linear combination,
with simplex coefficients, of point-like basis elements. (For $a=0$ we
have of course $\aacount{\\a0L} =\aacount{\\00L}
=\aacount{\\a0}\aacount{\\00L}$.)

\aaa{}
The remainder of this section provides a process to determine the
general $\aacount{\\abL}$, from the point-like elements of the form
$\aacount{\\00L^{(i)}}$. We have already done $b=0$. For $b\geq1$ we
need a new method. We look at products.

\aaa{}
We start at $b=1$. Consider the product
$\aacount{\\10L}\aacount{\\a0M}$, with $a\geq1$. Using
(\ref{eqn:[a0L]}) this product unpacks into:
\begin{equation}
  \label{eqn:[10L][a0M]}
  (
  \taacount{\\10}\taacount{\\00L}
  -
  \taacount{\\00}\taacount{\\00L^{(1)}}
  )
  (
  \taacount{\\a0}\taacount{\\00M}
  -
  \taacount{\\{a-1}0}\taacount{\\00M^{(1)}}
  )  
\end{equation}

\aaa{}
In (\ref{eqn:[10L][a0M]}) expanded, the coefficient of
$\aacount{\\00L}\aacount{\\00M}$ is $\aacount{\\10}\aacount{\\a0} =
\aacount{\\{a+1}0} + \aacount{\\{a-1}1}$. All other rank-zero
coefficients are of the form $\aacount{\\n0}$. Thus, we get an
equation with one `unknown'.  To simplify the calculations, we first
study the coefficients, in isolation.

\aaa{}
Taking a slightly more general starting point, we have
\begin{equation}
  \label{eqn:([a0]-[{a-1}0])-([b0]-[{b-1}0])}
  \begin{aligned}
    (\taacount{\\a0} - \taacount{\\{a-1}0})
    (\taacount{\\b0} - \taacount{\\{b-1}0})
    &=
    \taacount{\\{a+b}0}
    -2\taacount{\\{a+b-1}0}
    \taacount{\\{a+b-2}0}
    \\
    &\quad
    + \taacount{\\01}
    (\taacount{\\{a-1}0} - \taacount{\\{a-2}0})
    (\taacount{\\{b-1}0} - \taacount{\\{b-2}0})
  \end{aligned}
\end{equation}
assuming $a,b\geq0$.

\aaa{}
Now write
\begin{equation}
  \taacount{\\a0{\sim}} = \taacount{\\a0} - \taacount{\\{a-1}0}
\end{equation}
for $a\geq1$. This is similar to (\ref{eqn:[a0L]}). We call
$\taacount{\\a0{\sim}}$ a \emph{quasi-simplex}. It is part of the
\emph{shadow basis} for $\calR_{[0]}$. Note that it has mixed degree.

\aaa{}
Equation (\ref{eqn:([a0]-[{a-1}0])-([b0]-[{b-1}0])}) now becomes
\begin{equation}
  \taacount{\\a0{\sim}}\taacount{\\b0{\sim}}
  =
  \taacount{\\{a+b}0{\sim}}
  - \taacount{\\{a+b-1}0{\sim}}
  + D\taacount{\\{a-1}0{\sim}}\taacount{\\{b-1}0{\sim}}
\end{equation}
which holds for $a,b\geq 1$. This is the \emph{quasi-simplex product
  formula}. We use $D$ rather than $\aacount{\\01}$ to emphasise the
similarity with the cone product formula. [Exercise]

\aaa{}
Now consider the special case $a=b=1$. Our new product formula reduces
to:
\begin{equation}
  \label{eqn:[10sim][10sim]}
  \taacount{\\10{\sim}}\taacount{\\10{\sim}}
  =
  \taacount{\\20{\sim}}
  - \taacount{\\10{\sim}}
  + D
\end{equation}
Part of our process is to extend the shadow basis to a counting basis
for $\calR_{[0]}$. Bearing (\ref{eqn:[10sim][10sim]}) in mind, how
should we express $D$ in the shadow basis?

\aaa{}
How to write $D$ in the shadow basis is the unknown.  Because $D$ is
not a linear combination of simplices, nor is it of
quasi-simplices. So it must use a new basis vector, which we will call
$\taacount{\\01{\sim}}$. In addition, (\ref{eqn:[10sim][10sim]})
contains a negative. This must be cancelled, and only $D$ is available
to do this.

\aaa{}
This leads to the equation
\begin{equation}
  D = \taacount{\\01} = \taacount{\\01{\sim}} + \taacount{\\10{\sim}}
\end{equation}
which defines $\aacount{\\01{\sim}}$ in the shadow basis. It also leads to
\begin{equation}
  \taacount{\\01}\taacount{\\00L}
  =
  \taacount{\\01L} + \taacount{\\10L^{(1)}}
\end{equation}
which determines $\aacount{\\01L}$ in the candidate counting basis
(given the point-like elements).

\aaa{}
When applied to $\aacount{\\10{\sim}}\aacount{\\{n+1}0{\sim}}$ this
process works in exactly the same way, to give
\begin{align}
  \taacount{\\01}\taacount{\\n0{\sim}}
  &=
  \taacount{\\n1{\sim}} + \taacount{\\{n+1}0{\sim}}
  \\
  \taacount{\\01}\taacount{\\n0L}
  &=
  \taacount{\\n1L} + \taacount{\\{n+1}0L^{(1)}}
\end{align}
which respectively define $\aacount{\\n1{\sim}}$ and
$\aacount{\\n1L}$. [Exercise]

\aaa{}
Now consider:
\begin{align}
  \nonumber
  \taasim20 \taasim20
  &=
  \taasim40 - \taasim30 + D\taasim10\taasim10
  \\
  \label{eqn:[20sim][20sim]}
  &=
  \taasim40 - \taasim30 + D\taasim20 - D\taasim10 + D^2
\end{align}

\aaa{}
The big unknown is the shadow basis.  The unknown here is $D^2$ (in
the shadow basis). The first negative (\ref{eqn:[20sim][20sim]}) is
$-\aasim30$. However, the next term is $D\aasim20$, which is equal to
$\aasim21 + \aasim30$. Thus, the first negative is cancelled. The
remaining negative is $-D\aasim10 = -\aasim11 - \aasim20$.

\aaa{}
Thus, using the previous logic, we have
\begin{align}
  D^2
  &=
  \taasim02 + \taasim11 + \taasim20
  \\
  \taacount{\\02}\taacount{\\00L}
  &=
  \taacount{\\02L} + \taacount{\\11L^{(1)}} + \taacount{\\20L^{(2)}}
  \\
  \taacount{\\02}\taacount{\\n0L}
  &=
  \taacount{\\n2L} + \taacount{\\{n+1}1L^{(1)}} + \taacount{\\{2+n}0L^{(2)}}
\end{align}
where the last equation comes from $\aasim20 \aasim{2+n}0$.  This
process can be continued to provide a definition for all
$\aacount{\\abL}$, by using recursive on $b$. [Exercise].

\aaa{}
For clarity, this construction of the $\aacount{\\abL}$ depends on
knowing the point-like elements. This we do in the next section.

\section{Cone shadows}
\label{sec:cone-shadows}

\aaa{}
\label{par:[00abL]}
Here is the construction for the point-like
$\aacount{\\00L}$, for $r=\rank L\geq 1$.  Assume $L=L_1L'$, and
$L=\aatop{a}{b}$. The formula is
\begin{equation}
  \taacount{\\00\\abL}
  =
  C_{[1]}(\taacount{\\{a}{b+1}L})
  - C_{[1]}(\taacount{\\{a+2}{b}L})
\end{equation}
where the extra term $C_{[1]}(\taacount{\\{a+2}{b} L})$ can give rise
to \emph{shadow terms} in $C_{[1]}(\aacount{L})$.

\aaa{}
The extra term is present to ensure consistency with the calculations
in \cite[Prop.~19]{jf:g-vector-axioms}. If $b=0$ then
$\aacount{\\{a+2}{b}L}$ is simplex-like, and so the extra term is
zero.

\aaa{}
The understanding and explanation of this is outside the scope of the
present paper. The author expects this feature to be necessary for the
basis we have just constructed to be a counting basis. See
\S\ref{sec:advice} for some ideas as how that might be proved.


\begin{thebibliography}{10}

\bibitem
  {bb:gds}
  Margaret Bayer and Louis~J. Billera.
  Generalized {D}ehn-{S}ommerville relations for polytopes, spheres and
  {E}ulerian partially ordered sets.
  {\em Invent. Math.}, 79(1):143--157, 1985.

\bibitem%
  {jf:MV-IC}
  Jonathan Fine.
  The {M}ayer-{V}ietoris and {$IC$} equations for convex polytopes.
  {\em Discrete Comput. Geom.}, 13(2):177--188, 1995.

\bibitem%
  {jf:g-vector-axioms} Jonathan Fine,
  Axioms for the $g$-vector of general convex polytopes,
  arXiv:1011.4269 

\bibitem
  {jf:cpr} Jonathan Fine, Cone product rings, (in
  preparation)

\bibitem
  {jf:cplh} Jonathan Fine, Convex polytopes and linear homology, (in
  preparation)

\bibitem%
  {pm:num-faces} Peter McMullen, The numbers of faces of simplicial
  polytopes, {\em Israel J. Math.}, 9:559-570, 1971.

\bibitem
  {hk:betti-produkt} Hermann K\"unneth, \"Uber die Bettischen Zahlen
  einer Produktmannigfaltigkeit, {\em Mathematische Annalen},
  90:65-85, 1923.
  
\end{thebibliography}
\end{document}